\documentclass[12pt]{article}

\usepackage{amsfonts}
\usepackage{amsmath}
\usepackage{amssymb}
\usepackage{color}
\usepackage{amsthm}
\usepackage{hyperref}

\newtheorem{theorem}{Theorem}

\newtheorem{example}[theorem]{Example}

\newtheorem{proposition}[theorem]{Proposition}
\newtheorem{remark}[theorem]{Remark}

\newcommand{\R}{\mathbb R}
\newcommand{\N}{\mathbb N}
\numberwithin{theorem}{section}
\numberwithin{equation}{section}

\title{Simple variance bounds with applications to Bayesian posteriors and intractable distributions} \author{Fraser
  Daly\footnote{Department of Actuarial Mathematics and Statistics,
    Heriot-Watt University, Edinburgh EH14 4AS, UK.  E-mail:
    f.daly@hw.ac.uk},\;
  Fatemeh Ghaderinezhad\footnote{Department of Applied Mathematics,
    Computer Science and Statistics, Ghent University, Krijgslaan 281,
    S9, Campus Sterre, 9000 Gent, Belgium. E-mail:
    fatemeh.ghaderinezhad@ugent.be}, Christophe
  Ley\footnote{Department of Applied Mathematics, Computer Science and
    Statistics, Ghent University, Krijgslaan 281, S9, Campus Sterre,
    9000 Gent, Belgium. E-mail: christophe.ley@ugent.be}\; and Yvik
  Swan\footnote{Department of Mathematics, Universit\'e libre de
    Bruxelles, Boulevard du Triomphe, CP210, B-1050 Bruxelles}} \date{\today}

\begin{document}

\maketitle

\noindent{\bf Abstract} Using coupling techniques based on Stein's method for probability approximation, we revisit classical variance bounding inequalities of Chernoff, Cacoullos, Chen and Klaassen. Taking advantage of modern coupling techniques allows us to establish novel variance bounds in settings where the underlying density function is unknown or intractable. Applications include bounds for asymptotically Gaussian random variables using zero-biased couplings, bounds for random variables which are New Better (Worse) than Used in Expectation, and analysis of the posterior in Bayesian statistics.   \vspace{12pt}


\vspace{12pt}


\section{Introduction}
\label{sec:literature-review}

Weighted Poincar\'e (or isoperimetric) inequalities, giving upper
bounds on the variance of a function of a random variable, have a long
and rich history, beginning with the work of Chernoff
\cite{Chernoff81}. Chernoff proved that if $X$ has a centred Gaussian
distribution with variance $\sigma^2$, then
\begin{equation}
  \label{eq:14}
  \mathrm{Var}[g(X)]\leq \sigma^2\mathbb{E}[(g^\prime(X))^2]\,,
\end{equation}
for any absolutely continuous function $g:\mathbb{R}\mapsto\mathbb{R}$
such that $g(X)$ has finite variance. This inequality has since been
generalized by many authors, including Cacoullos \cite{C82}, Chen
\cite{Chen82} and Klaassen \cite{Klaassen85}. To accompany these upper
variance bounds, many of these authors have also established
corresponding lower bounds, in the form of generalized Cram\'er-Rao
inequalities. In particular in the centred Gaussian case we have
\begin{equation}
  \label{eq:12}
  \mathrm{Var}[g(X)]\geq \sigma^2\mathbb{E}[g^\prime(X)]^2,
\end{equation}
see \cite{C82}.  The above cited works represent early entries in what
is now a vast literature; we refer to \cite{ERS19vb1,ERSvb2} for
recent overviews of this large body of work.

The purpose of the present article is to revisit these classical
variance bounding inequalities in light of the coupling techniques at
the heart of Stein's method for probability approximation (see, for
example, \cite{ChGoSh11} and \cite{ley2017stein} for recent introductions to Stein's
method). These techniques allow us to establish upper and lower
variance bounds in a variety of settings, including many in which the
density of the underlying random variable is unknown or
intractable. Making use, for example, of the zero-biased coupling
allows us to establish explicit variance bounds for a wide range of
situations in which the underlying random variable is known to be
asymptotically Gaussian. In Sections
\ref{sec:appl-stein-fram}--\ref{sec:bounds-using-stoch} we will
consider a variety of situations where bounds may be derived using
this, and other, couplings. Before doing so, we use the remainder of
this section to outline the general coupling techniques we employ from
Stein's method, and how these can be used to establish upper and lower
variance bounds in the spirit of Chernoff, Cacoullos, Chen and
Klaassen.

 Let $W$ be a real random variable on some
  fixed probability space. Let $\gamma$ be a real-valued function.
 We say that a pair of random variables
  $(T_1, T_2)$ (living on the same probability space as $W$) form a
  Stein coupling for $W$ with respect to $\gamma$ if 
  \begin{equation}
    \label{eq:1}
    \mathbb{E} \left[ \gamma(W) \phi(W) \right] = \mathbb{E} \left[ T_1 \phi'(T_2) \right]
  \end{equation}
  for all test functions $\phi \in C$ with $C\subset C^1(\R)$ some
  appropriately chosen class of functions. Although the choice
  $C=C_0^{\infty}(\R)$ is always allowed, it will generally be
  necessary to use $C$ as wide as possible; this fact is often
  reflected in the literature wherein one rather makes use of the
  generic expression ``where $C$ is the class of functions for which
  expectations on both sides
  exist''.  

  We begin by showing an elementary argument allowing us to use
  \eqref{eq:1} to obtain tight upper variance bounds. To this end,
  suppose that $\gamma$ is a strictly
  increasing, differentiable function with exactly one sign change. Then in particular
  it is invertible and $\gamma^{-1}(0)$ is well-defined. Let $g$ be a
  real-valued differentiable function such that $\mathrm{Var}[g(W)]$ is
  finite. Following   \cite{LS16} we write 
  \begin{align*}
    \mathrm{Var}[g(W)] & \le \mathbb{E} \left[ \left( g(W) - g(\gamma^{-1}(0))
                       \right)^2 \right] = \mathbb{E} \left[ \left( \int_0^{\gamma(W)}
                       \frac{g'(\gamma^{-1}(u))}{\gamma'(\gamma^{-1}(u))}\mathrm{d}u \right)^2
                       \right] \\
                     & \le \mathbb{E} \left[ \gamma(W)  \int_0^{\gamma(W)}
                       \left( \frac{g'(\gamma^{-1}(u))}{\gamma'(\gamma^{-1}(u))}
                       \right)^2\mathrm{d}u \right]\,, 
  \end{align*}
  where the  equality follows by differentiability of $g$ and the
  subsequent inequality via Cauchy-Schwarz. Applying \eqref{eq:1} as well as
  Leibnitz' rule for differentiating integrals we deduce the general
  upper variance bound
  \begin{equation}
    \label{eq:2}
 \mathrm{Var}[g(W)]   \le \mathbb{E} \left[ \frac{T_1}{\gamma'(T_2)} \left(
     g'(T_2) \right)^2 \right], 
  \end{equation}
  which holds as soon as the function
  $x \mapsto \int_0^{\gamma(x)} \left(
    \frac{g'(\gamma^{-1}(u))}{\gamma'(\gamma^{-1}(u))} \right)^2\mathrm{d}u$
  belongs to the (so far unspecified) class $C$.  Note that inequality
  \eqref{eq:2} also holds if in \eqref{eq:1} we replace the equality
  sign by an increasing inequality.

  Identity \eqref{eq:1} can also readily be combined with the
  Cauchy-Schwarz inequality to obtain lower variance bounds. To this
  end, consider a mean zero function $\gamma$ (this is in any case
  necessary for relationships such as \eqref{eq:1} to hold) for which
  $(\mathbb{E} \left[ \gamma(W) g(W) \right])^2 = (\mathbb{E} \left[
    \gamma(W)(g(W) - \mathbb{E}[g(W)]) \right])^2 \le \mathbb{E}
  \left[ \gamma(W)^2 \right] \mathrm{{Var}}[g(W)]$.  Then from
  \eqref{eq:1} we deduce
  \begin{equation}
    \label{eq:3}
   \mathrm{Var}[g(W)]\ge \frac{(\mathbb{E} \left[ T_1 g'(T_2) \right])^2}{\mathrm{Var}
      \left[ \gamma(W) \right]}
  \end{equation}
  for all $g \in C$.  As above, we note that inequality \eqref{eq:3}
  also holds if in \eqref{eq:1} we replace the equality sign by a
  decreasing inequality.  
  
The rest of this paper is devoted to proposing situations wherein such couplings $W, T_1$ and $T_2$ occur naturally and may be used to establish upper and lower variance bounds. In Section \ref{sec:appl-stein-fram} we use the framework of Stein kernels to express suitable couplings. Section \ref{sec:vb-via-biasing} makes use of zero-biased couplings to derive variance bounds suitable for random variables which are asymptotically Gaussian. Finally, in Section \ref{sec:bounds-using-stoch} we consider random variables satisfying certain stochastic or convex ordering assumptions, which allow us to derive bounds sharper than we would otherwise obtain with our method. Some proofs and additional examples illustrating the results of Section \ref{sec:appl-stein-fram} are deferred to the appendices.

  \section{Stein kernel and a bound of Cacoullos}
\label{sec:appl-stein-fram}

Suppose that the target $W$ has a differentiable density $p$ with interval support.  Following, for example, \cite{CPU94} and \cite{ERS19vb1}, we define the \emph{Stein kernel} of $W$ as the function $\tau$ satisfying
\begin{equation}
  \label{eq:6}
\mathrm{Cov}\left[W,  \phi (W)\right] =  \mathbb{E} \left[
 \tau(W) \phi'(W) \right] 
\end{equation}
for all functions $\phi$ such that either integral is defined. See \cite{ERS19vb1} for an extensive discussion of this function.  In the notation
of Section \ref{sec:literature-review}, this means that we can take
$\gamma(x) = x-\mathbb{E}[W]$, $T_1 = \tau(W)$ and $T_2 = W$ in \eqref{eq:1}. Note
that $\mathbb{E}[\tau(W)] = \mathrm{Var}[W]$.  Applying \eqref{eq:2}
and \eqref{eq:3}, we get  for all $g \in L^2(W)$ that
  \begin{equation}
    \label{eq:8}
   \frac{ \mathbb{E} \left[ 
        \tau(W) g'(W)
       \right]^2}{\mathrm{Var} \left[  W  
       \right]} \le  \mathrm{{Var}}[g(W)]
     \le \mathbb{E} \left[ \tau(W) \left( g'(W) \right)^2
    \right], 
  \end{equation}
  which is nothing but a restatement of classical bounds already
  available in \cite{C82}.

  Of course for \eqref{eq:8} to be of use it remains to identify
  situations in which the Stein kernel has an agreeable
  form. We give several such situations.

  \begin{example}
    Following \cite{nourdin2013integration}, it is easy to see that if
    $W = n^{-1/2} \sum_{i=1}^n X_i$, where the $X_i$ are centred,
    independent random variables with Stein kernel $\tau_i(\cdot)$ and common variance $\sigma^2$,
    then
    $ \tau_W(w) = \frac{1}{n} \sum_{i=1}^n\mathbb{E}[\tau_i(X_i) \, |
    \, W = w]$ is a Stein kernel for $W$. If the $X_i$ are copies of
    $X_1$ with kernel $\tau_{1}(\cdot)$, \eqref{eq:8} becomes
  \begin{equation*}
    \frac{\mathbb{E}\left[ \tau_1(X_1) g'(W)\right]^2}{\sigma^2} \le
    \mathrm{Var}[g(W)] \le 
    \mathbb{E}\left[\tau_1(X_1) (g'(W))^2\right]
  \end{equation*}
  If $W$ and $X_1$ were independent, we could  use
  $\mathbb{E}[\tau_1(X_1)]= \sigma^2$ to recover  the
  Gaussian case stated in \eqref{eq:14} and
  \eqref{eq:12}. Here we need to apply a limited development to make 
independence appear.   Let $U \sim \mathrm{Unif}[0,1]$ and recall the mean-value theorem
  $g'(x+t) = g'(x) + t \mathbb{E}[g''(x+Ut)]$. Let
  $W^{(1)} = W - n^{-1/2}X_1$. Then, by independence, if $g$ is twice
  differentiable 
  the lower bound becomes
  $ \sigma^2\mathbb{E}[g'(W^{(1)})]^2 + \frac{C_1 }{\sqrt n}$ where
  $C_1 =C_1(g,n)$ is given by $C_1= 2\mathbb{E}[g'(W^{(1)})]\mathbb{E}[ \tau_1(X_1)
  X_1g''(W^{(1)}+n^{-1/2}UX_1)] + n^{-1/2}/\sigma^2\mathbb{E}[ \tau_1(X_1)
  X_1g''(W^{(1)}+n^{-1/2}UX_1)]^2$.  Clearly $\lim_{n\rightarrow\infty} C_1(g,n)/\sqrt{n}=0$ for all $g$. Similar considerations apply for
  the upper bound. Indeed, recall that for a twice differentiable function $g$
  we have
\begin{equation}
  \label{eq:WstarTVI}
  \left|g^\prime(x+t)^2-g^\prime(x)^2\right|\leq2\|
g^\prime 
g^{\prime\prime} \| |t|\,, 
\end{equation}
(where $\lVert\cdot\rVert$ is the supremum norm) so that we have
$\mathbb{E}[\tau_1(X_1) (g'(W))^2] \le \sigma^2
\mathbb{E}[(g'(W^{(1)}))^2] + \frac{2}{\sqrt n} \| g' g''\| \sigma^2
\mathbb{E}[|X_1|] =: \sigma^2 \mathbb{E}[(g'(W^{(1)}))^2] +
\frac{C_2}{\sqrt n}$.  Wrapping up,
  \begin{equation*}
 \sigma^2
    \mathbb{E}[(g'(W^{(1)}))]^{2} + \frac{C_1}{\sqrt n} \le     \mathrm{Var}[g(W)] \le \sigma^2
    \mathbb{E}[(g'(W^{(1)}))^2] + \frac{C_2}{\sqrt n},
  \end{equation*}
  where the proximity with the corresponding inequalities for the
  Gaussian case are now made explicit.
\end{example}

\begin{example}[Smoothing] \label{eg:smoothing} Let $Y$ be a
  real-valued random variable with $\mathbb{E}[Y]=\mu$. Note that we
  do not require $Y$ to have a density function, and the bounds of
  this example apply if, for instance, $Y$ is a discrete random
  variable. In order to allow us to derive variance bounds for $Y$
  using our approach, we smooth it by convolving it with independent
  Gaussian noise with small variance. We let
  $Z\sim\mathcal{N}(0,\epsilon^2)$ have a Gaussian distribution,
  independent of $Y$.  Let $\varphi_\epsilon$ and $\Phi_\epsilon$ be
  the density and distribution functions of $Z$, respectively, and
  define
\begin{equation}\label{eq:smooth}
  \tau_{\epsilon}(x)=\epsilon^2+\frac{\mathbb{E}\left[(Y^\prime-\mu)\bar{\Phi}_\epsilon(x-Y^\prime)\right]}{\mathbb{E}[\varphi_\epsilon(x-Y^\prime)]}\,, 
\end{equation}
where $\bar{\Phi}_\epsilon(y)=1-\Phi_\epsilon(y)$ and $Y^\prime$ is an
independent copy of $Y$. Then $\tau_{\epsilon}(x)$ is a Stein kernel
for $Y+Z$ (see Appendix \ref{sec:proofs}) and \eqref{eq:8} applies to
all differentiable functions $g:\mathbb{R}\mapsto\mathbb{R}$ such that
$\mathrm{Var}[g(Y+Z)]$ is finite. Moreover, the following hold:
\begin{itemize}
\item[(i).] If the mapping
  $x\mapsto\left(g(x)-\mathbb{E}[g(Y+Z)]\right)^2$ is convex, then 
$$
\mathrm{Var}[g(Y)]\leq\mathbb{E}\left[\tau_{\epsilon}(Y+Z)g^\prime(Y+Z)^2\right]\,.
$$
\item[(ii).] If the mapping $x\mapsto (g(x)-\mathbb{E}[g(Y)])^2$ is concave, then
$$
\mathrm{Var}\left[g(Y)\right]\geq\frac{\mathbb{E}\left[\tau_{\epsilon}(Y+Z)g^\prime(Y+Z)\right]^2}{\epsilon^2+\mathrm{Var}[Y]}\,.
$$
\end{itemize}
We defer the proofs of these claims to Appendix \ref{sec:proofs}.
\end{example}

\begin{example}[Pearson family and application to posterior distributions] \label{eg:Pearson}

  As is well known, the Pearson family has explicit Stein kernels
  given by Proposition \ref{prop:perason} recalled in the Appendix.
  Such a result is particularly useful in the following situation
inherited from Bayesian statistics.  In a Bayesian setting, the
  initial distribution of the parameter of interest is some prior distribution with density
  $ \pi_0(\theta)$; upon observing data points
  $\mathbf{x}=(x_1, \ldots, x_n)$ sampled independently with sampling
  distribution $\pi(\theta, \mathbf{x})$ we update from the prior to
  the posterior density given by
  $\pi_2(\theta) = \kappa_2(\mathbf{x})\pi(\theta,
  \mathbf{x})\pi_0(\theta)$. We use the notation $\Theta_0$ to
  indicate the distribution of the parameter under the prior,
  $\Theta_2$ its distribution under the posterior, and $X$ a random
  variable following the same common distribution of the
  observations. We also write $\Theta_1$ for the parameter under the sampling
  distribution
  $\pi_1(\theta) = \kappa_1(\mathbf{x}) \pi(\theta, \mathbf{x})$, which
  corresponds to a posterior with flat (uninformative) prior. A
  popular choice of prior is that of a \emph{conjugate} prior for
  which the mathematical properties of the posterior are the same as
  those of the sampling distribution; the impact of the data is then
  visible in the parameters of the posterior distribution who are updated. Restricting
  our attention to Pearson distributed families, we can apply
  Proposition \ref{prop:perason} and read variance bounds directly
  from the updated parameters. For instance:
  \begin{itemize}
  \item Gaussian data, inference on mean, Gaussian prior:    If $X \sim \mathcal{N}(\theta, \sigma^2)$ with $\theta\in\R$ and
  fixed $\sigma>0,$ and $\Theta_0 \sim \mathcal{N}(\mu, \delta^2)$
  with $\mu\in\R,\delta>0$, then
  $\Theta_2 \sim
  \mathcal{N}\left(\frac{\sigma^2\mu+n\delta^2\bar{x}}{n\delta^2+\sigma^2},\frac{\sigma^2\delta^2}{n\delta^2+\sigma^2}\right)$,
  where $\bar{x}=\frac{1}{n}\sum_{i=1}^nx_i$. The Stein kernel for
  this Gaussian distribution is
  $\tau(\theta) = ( \frac{n}{\sigma^2} +
  \frac{1}{\delta^2})^{-1}$. Consequently, 
$$
\mathbb{E} \left[g'(\Theta_2) \right]^2 \leq \bigg( \frac{n}{\sigma^2} +
  \frac{1}{\delta^2}\bigg) \mathrm{Var}[g(\Theta_2)] \leq \mathbb{E} [g'(\Theta_2)^2]
  $$
  for all suitable $g$, all $n$ and all values of the parameters.
\item   Gaussian data, inference on variance, Inverse Gamma  prior:  If  $X \sim \mathcal{N}(\mu, \theta)$ with $\theta>0$ and fixed $\mu\in\R$, and $\Theta_0 \sim \mathcal{IG}(\alpha, \beta)$ has an Inverse Gamma distribution with density $$
\theta \mapsto \frac{\beta^\alpha}{\Gamma(\alpha)} \theta^{-\alpha -1} \exp \left(  -\frac{\beta}{\theta} \right), \, \alpha,\beta>0,
$$ then
   $\Theta_2 \sim \mathcal{IG}\left(\frac{n}{2} + \alpha, \frac{1}{2}
     \sum_{i=1}^n (x_i - \mu)^2 + \beta\right)$.  The Stein kernel for
   this Inverse Gamma distribution is   $\tau (\theta) =
   \frac{\theta^2}{\frac{n}{2}+\alpha-1}$. Consequently, for all suitable $g$, 
$$
\frac{(\frac{n}{2} + \alpha-2)}{(\frac{1}{2} \sum_{i=1}^n (x_i - \mu)^2 + \beta)^2} \mathbb{E} [\Theta_2^2g'(\Theta_2)]^2 \leq \mathrm{Var}[g(\Theta_2)] \leq \frac{1}{\frac{n}{2}+\alpha-1} \mathbb{E} [\Theta_2^2g'(\Theta_2)^2].
$$

\item Binomial data, inference on proportion, Beta  prior: If  $X \sim {Bin}(n, \theta)$ with $\theta\in[0,1]$, and $\Theta_0 \sim{Beta}(\alpha, \beta)$ with density $$
\theta \mapsto \frac{\theta^{\alpha-1}(1-\theta)^{\beta-1}}{\frac{\Gamma(\alpha)\Gamma(\beta)}{\Gamma(\alpha+\beta)}}, \, \alpha,\beta>0,
$$ then
$\Theta_2 \sim {Beta}\left(x+\alpha,n-x+\beta\right)$, where $x$
denotes the observed number of successes.  The Stein kernel for this
Beta distribution is
$\tau (\theta) = \frac{\theta (1 - \theta)}{n+ \alpha +
  \beta}$. Consequently, for all suitable $g$,
$$
\frac{(n+\alpha +\beta+1)}{(x+\alpha )(n - x + \beta)} \mathbb{E} [\Theta_2 (1-\Theta_2) g'(\Theta_2)]^2 \leq \mathrm{Var}[g (\Theta_2)] \leq \frac{\mathbb{E}[\Theta_2 (1 - \Theta_2)g'(\Theta_2)^2]}{n+\alpha +\beta}.
$$ 
\end{itemize}

Further examples are provided in Appendix \ref{sec:more-examples-1}. 

\end{example}

\section{Variance bounds from zero-biased couplings}
\label{sec:vb-via-biasing}

In this section, we suppose that the target $W$ has mean zero, finite variance
$\sigma^2$, and  can be coupled to some
random variable $W^{\star}$ through
\begin{equation}
\label{eq:bias}
  \mathbb{E}[W\phi(W)] = \sigma^2 \mathbb{E}   [\phi'(W^\star)]
\end{equation}
for all functions $\phi:\mathbb{R}\mapsto\mathbb{R}$.  Such $W^{\star}$ always exists, and its law is unique. It
has the $W$-zero-biased distribution; see, e.g., \cite[Section
2.3.3]{ChGoSh11} and references therein for more details.  Note that
$W^\star$ is a continuous random variable, regardless of whether $W$
is discrete or continuous.
Under \eqref{eq:bias}, we
immediately obtain 
\begin{equation}\label{eq:zb1}
\sigma^2\mathbb{E}\left[g^\prime(W^\star)\right]^2\le
\mathrm{Var}[g(W)]\leq\sigma^2\mathbb{E}\left[g^\prime(W^\star)^2\right]
\end{equation}
by using \eqref{eq:2} and \eqref{eq:3} with $\gamma(x) = x$,
$T_1 = \sigma^2$ and $T_2 = W^{\star}$ for all
$g:\mathbb{R}\mapsto\mathbb{R}$ for which $\mathrm{Var}[g(W)]$ is
finite.  Obviously it may be of interest to express \eqref{eq:zb1} in
terms of the original variable. Using \eqref{eq:WstarTVI}, we obtain
the following result.
\begin{proposition}
  Let $W$ have mean zero and finite variance $\sigma^2$, and $W^{\star}$ have the
  $W$-zero biased distribution.  Then
  \begin{equation}\label{eq:bdd}
  \mathrm{Var}[g(W)]\leq\sigma^2\mathbb{E}\left[g^\prime(W)^2\right]+
  2\sigma^2\|g^\prime
  g^{\prime\prime} \|\mathbb{E}|W^\star-W|
\end{equation}
for all twice differentiable functions $g:\mathbb{R}\mapsto\mathbb{R}$
for which $\mathrm{Var}[g(W)]$ exists. 
\end{proposition}

It is classical that the Gaussian distribution is the unique fixed
point of the zero-bias transform, in the sense that
$W\sim \mathcal{N}(0, \sigma^2)$ if and only if $W = W^{\star}$. Hence
$|W^\star - W|$  gives
information on the distributional proximity between the law $\mathcal{L}(W)$ of $W$
and $\mathcal{N}(0, \sigma^2)$. Also,
it is classical that the Gaussian is characterized by the fact that
$  \sigma^2 = \sup_g {\mathrm{Var}[g(W)]}/{\mathbb{E}[g'(W)^2]}$, see, e.g., \cite{CPU94}. Inequality \eqref{eq:bdd} captures these two
essential 
features of the Gaussian distribution.

\begin{example}
  Let $X_1,X_2,\ldots,X_n$ be 
  independent mean zero  random variables with finite variances
  $\mathbb{E}[X_i^2] = \sigma^2_i, i=1,\ldots,n$. Set $W = X_1 + \cdots + X_n$ and
  $\mathbb{E}[W^2] = \sigma^2=\sum_{i=1}^n\sigma_i^2$.  Let $I$ be a random index independent
  of all else such that $P(I = i) = \sigma^2_i/\sigma^2$ and let $W_i
  = W-X_i$. Finally let $X_i^{\star}$ be the zero-bias transform of
  $X_i$. Then  $ W^{\star}-W = X_I - X_I^{\star}$ (see Example 2.1 of \cite{GR97}) so that the bound
  \eqref{eq:bdd} becomes
  \begin{align*}
    \mathrm{Var}[g(W)] 
    & \le \sigma^2 \mathbb{E}[g'(W)^2]
                         + 2 \|g'g''\| \sum_{i=1}^n \sigma^2_i \mathbb{E}[|X_i - X_i^{\star}|] . 
  \end{align*}
  If, furthermore, we suppose the summands to be independent copies of
  $X$ such that  $\sigma^2 = 1$ 
  then
  \begin{align*}
    \mathrm{Var}[g(W)]  \le \mathbb{E}[g'(W)^2]
    + 2 \|g'g''\|   \mathbb{E}[|X - X^{\star}|]\,.
  \end{align*}
  To see how this plays out in practice, suppose that
  $X = (\xi-p)/\sqrt{npq}$ with $\xi$ Bernoulli with success parameter
  $p$. Following \cite[Corollary 4.1]{ChGoSh11}, we obtain
  $\mathbb{E}[|X - X^{\star}|] = (p^2+q^2)/(2\sqrt{npq})$ and 
\begin{equation*}
  \mathrm{Var}[g(W)] \le \sigma^2 \mathbb{E}[g'(W)^2]
    +   \|g'g''\|  \frac{p^2+q^2}{\sqrt{npq}}.
\end{equation*}
Many other examples can be explicitly worked out along these lines. 
\end{example}

\begin{example}
Let $(a_{i,j})_{i,j=1}^n$ be an array of real numbers and $\pi$ a uniformly chosen permutation of $\{1,\ldots,n\}$.  Let $W=\sum_{i=1}^na_{i,\pi(i)}$.  We further define
$$
a_{\bullet\bullet}=\frac{1}{n^2}\sum_{i,j=1}^na_{i,j}\,,\quad
a_{i\bullet}=\frac{1}{n}\sum_{j=1}^na_{i,j}\,,\quad\mbox{and}\quad
a_{\bullet j}=\frac{1}{n}\sum_{i=1}^na_{i,j}\,,
$$
and note that $\mathbb{E}[W]=na_{\bullet\bullet}$ and
$$
\mathrm{Var}[W]=\sigma^2=\frac{1}{n-1}\sum_{i,j=1}^n\left(a_{i,j}-a_{i\bullet}-a_{\bullet j}+a_{\bullet\bullet}\right)^2\,.
$$
See, for example, \cite[Section 4.4]{ChGoSh11}.  Letting
$Z=\sigma^{-1}(W-na_{\bullet\bullet})$ and
$C=\max_{1\leq i,j\leq n}|a_{i,j}-a_{i\bullet}-a_{\bullet
  j}+a_{\bullet\bullet}|$, the proof of Theorem 6.1 of \cite{ChGoSh11}
shows that $\mathbb{E}|Z^\star-Z|\leq8C\sigma^{-1}$ for some positive constant $C$, and so we have
from (\ref{eq:bdd}) that
$$
\mathrm{Var}[g(Z)]\leq\mathbb{E}\left[g^\prime(Z)^2\right]+\frac{16C}{\sigma}\|g^\prime
g^{\prime\prime}\|\,,
$$ 
for all twice differentiable $g$ such that $\mathrm{Var}[g(Z)]$ is finite.
\end{example}

\section{Variance bounds using stochastic ordering}
\label{sec:bounds-using-stoch}

We consider now some further applications in which we do not require explicit knowledge of the
density of $W$ in order to derive bounds on $\mathrm{Var}[g(W)]$ using our techniques. Unlike those examples in Section \ref{sec:vb-via-biasing}, the bounds we obtain here have the same form as in applications where we employ the exact expression for the underlying density, as in Section \ref{sec:appl-stein-fram}, without any additional `remainder' terms. We may obtain such bounds under natural assumptions on the random variable $W$, which we express in terms of stochastic orderings; the price we pay is in some restriction on the class of functions $g$ for which the bounds apply.

We begin by recalling the definitions of the orderings which we will use. For any random
variables $X$ and $Y$, we will say that $X$ is stochastically smaller
than $Y$ (denoted $X\leq_{st}Y$) if
$\mathbb{P}(X>t)\leq\mathbb{P}(Y>t)$ for all $t$.  We will say that
$X$ is smaller than $Y$ in the convex order (denoted $X\leq_{cx}Y$) if
$\mathbb{E}[\phi(X)]\leq\mathbb{E}[\phi(Y)]$ for all convex functions
$\phi$ for which the expectations exist. See \cite{ss07} for background and many
further details.

\subsection{Zero-biased couplings and the convex order}
\label{sec:ConvexOrder}

Let $W$ be a real-valued random variable with mean zero and variance $\sigma^2$. Recall the definition (\ref{eq:bias}) of $W^\star$, the zero-biased version of $W$. We note that, from Lemma 2.1(ii) of \cite{GR97}, $W^\star$ is supported on the closed convex hull of the support of $W$ and has density function given by 
\begin{equation}\label{eq:density}
p^\star_W(w)=\frac{1}{\sigma^2}\mathbb{E}[WI(W>w)]\,.
\end{equation}

If we assume that $W^\star\leq_{cx}W$, then we may write
\begin{equation}\label{eq:conv}
\mathbb{E}[W\phi(W)]=\sigma^2\mathbb{E}[\phi^\prime(W^\star)]\leq\sigma^2\mathbb{E}[\phi^\prime(W)]\,,
\end{equation}
for all differentiable functions $\phi$ such that $\phi^\prime$ is convex. That is, (\ref{eq:1}) holds with the equality replaced by an inequality for all such $\phi$, with the choices $\gamma(W)=W$, $T_1=\sigma^2$, and $T_2=W$. 

Following the proof of (\ref{eq:2}), the inequality (\ref{eq:conv}) is sufficient to obtain this upper bound on $\mathrm{Var}[g(W)]$. In proving this bound, we apply (\ref{eq:conv}) with $\phi$ such that $\phi^\prime(x)=g^\prime(x)^2$; we must therefore assume that $g^\prime(x)^2$ is convex in order to do this. We thus obtain the following bound.
\begin{theorem}\label{theorem:convex}
Let $W$ have mean 0 and variance $\sigma^2$, and assume that $W^\star\leq_{cx}W$. For all differentiable $g:\mathbb{R}\mapsto\mathbb{R}$ such that $\mathrm{Var}[g(W)]$ exists and $g^\prime(x)^2$ is convex,
\begin{equation}\label{eq:convex_bd}
\mathrm{Var}[g(W)]\leq\sigma^2\mathbb{E}[g^\prime(W)^2]\,.
\end{equation}
\end{theorem}    
\begin{example}
Let $W=X_1+X_2+\cdots+X_n$, where $X_1,X_2,\ldots,X_n$ are independent, mean-zero random variables, with $X_i$ supported on the set $\{-a_i,b_i\}$ for $a_i,b_i>0$, for each $i=1,\ldots,n$. That is, $\mathbb{P}(X_i=-a_i)=p_i=1-\mathbb{P}(X_i=b_i)$ for $1\leq i\leq n$, where $p_i=b_i/(a_i+b_i)$ so that $\mathbb{E}[X_i]=0$.  Let $\sigma_i^2=\mathrm{Var}(X_i)$ and $\sigma^2=\sigma_1^2+\cdots+\sigma_n^2$.  

A straightforward calculation using (\ref{eq:density}) shows that, for each $i=1,\ldots,n$, $X_i^\star$ is uniformly distributed on the interval $[-a_i,b_i]$. Hence, Theorem 3.A.44 of \cite{ss07} gives that $X_i^\star\leq_{cx}X_i$ for each $i$.

Let $I$ be a random index, chosen independently of all else, with $\mathbb{P}(I=i)=\sigma_i^2/\sigma^2$, for $i=1,\ldots,n$. Now, using Lemma 2.1(v) of \cite{GR97}, $W^\star$ is equal in distribution to $X_I^\star+\sum_{j\not=I}X_j$, which is smaller than $W$ in the convex order for each possible value of $I$ by (3.A.46) of \cite{ss07}. It then follows from Theorem 3.A.12(b) of \cite{ss07} that $W^\star\leq_{cx}W$, and hence our upper bound (\ref{eq:convex_bd}) applies.
\end{example}

\subsection{Equilibrium couplings}
\label{sec:StocOrder}

Throughout this section, let $W$ be a non-negative random variable with mean $\lambda^{-1}$. Following, for example, \cite{PeRo11}, we say that a random variable $W^e$ has the equilibrium distribution with
respect to $W$ if 
\begin{equation}\label{eq:EqCoup}
\mathbb{E}[\phi(W)]-\phi(0)=\lambda^{-1}\mathbb{E}[\phi^\prime(W^e)]\,,
\end{equation}
for all a.e.\ differentiable functions $\phi$.

\begin{remark}
  Note that this definition is motivated by the fact that $W$ is
 Exponential if and only if $W$ and $W^e$ are equal in distribution.
  Applying the definition to the function
  $\phi_x(w) = (w-x) \mathbb{I}(w \ge x)$ and integrating by parts we
  obtain that
  $\mathbb{P}(W^e>x)=\lambda\int_x^\infty\mathbb{P}(W>y)\,\mathrm{d}y$
  for all $x \ge 0$.
\end{remark}

In this section we consider random variables that are new better than used in expectation (NBUE) and new worse than used in expectation (NWUE). Recall that $W$ is NBUE if $\lambda\int_x^\infty \mathbb{P}(W>s)\,ds\leq\mathbb{P}(W>x)$ for all $x\geq0$, and that $W$ is NWUE if this holds with the inequality reversed. These properties are well-known in reliability theory; see, for example, \cite{ss07}.

From this definition and the remark above, it is clear that $W$ is NBUE if and only if $W^e\leq_{st}W$,  and that $W$ is NWUE if and only if $W\leq_{st}W^e$. For a random variable $W$ which is either NBUE or NWUE, we employ this stochastic ordering in a similar way to the convex ordering we used in Section \ref{sec:ConvexOrder} above.

We begin by deriving an inequality analogous to (\ref{eq:conv}). For a differentiable function $\phi$, the definition of $W^e$ gives that
$$
\mathbb{E}[W\phi(W)]=\lambda^{-1}\mathbb{E}[\phi(W^e)+W^e\phi^\prime(W^e)]\,,
$$
and hence
$$
\mathbb{E}[(\lambda W-1)\phi(W)]+\mathbb{E}[\phi(W)]=\mathbb{E}[W^e\phi^\prime(W^e)]+\mathbb{E}[\phi(W^e)]\,.
$$
Thus, the inequality
\begin{equation}\label{eq:equilibrium}
\mathbb{E}[(\lambda W-1)\phi(W)]\leq\mathbb{E}[W\phi^\prime(W)]
\end{equation}
holds if and only if
$$
\mathbb{E}[\phi(W^e)+W^e\phi^\prime(W^e)]\leq\mathbb{E}[\phi(W)+W\phi^\prime(W)]\,.
$$
Therefore, inequality (\ref{eq:equilibrium}) holds if $W$ is NBUE and $\phi(x)+x\phi^\prime(x)$ is increasing in $x$. Alternatively, (\ref{eq:equilibrium}) also holds if $W$ is NWUE and $\phi(x)+x\phi^\prime(x)$ is decreasing in $x$. Analogously to the use of (\ref{eq:conv}) in proving Theorem \ref{theorem:convex} above, an upper bound on $\mathrm{Var}[g(W)]$ therefore holds for some functions $g$ under either of these assumptions; see Theorem \ref{theorem:equilibrium} below for a precise statement.

Similarly, we may ask when the reversed inequality $\mathbb{E}[(\lambda W-1)\phi(W)]\geq\mathbb{E}[W\phi^\prime(W)]$ holds. By similar reasoning, this holds if either (i) $W$ is NBUE and $\phi(x)+x\phi^\prime(x)$ is decreasing in $x$, or (ii) $W$ is NWUE and $\phi(x)+x\phi^\prime(x)$ is increasing in $x$. Under either of these assumptions, we have a lower variance bound.

We have thus proved the following.
\begin{theorem}\label{theorem:equilibrium}
Let $W$ be a non-negative random variable with mean $\mathbb{E}[W]=\lambda^{-1}$.
\begin{enumerate}
\item[(a)] For a differentiable function $g:\mathbb{R}^+\mapsto\mathbb{R}$ such that $\mathrm{Var}[g(W)]$ exists, let $\phi_g(x)=\int_0^{\lambda x-1}g^\prime(\lambda^{-1}(u+1))\,du$. Assume that either
\begin{enumerate}
\item[(i)] $W$ is NBUE and $\phi_g(x)+x\phi_g^\prime(x)$ is increasing in $x$; or
\item[(ii)] $W$ is NWUE and $\phi_g(x)+x\phi_g^\prime(x)$ is decreasing in $x$.
\end{enumerate}
Then
$$
\mathrm{Var}[g(W)]\leq\frac{1}{\lambda}\mathbb{E}[Wg^\prime(W)^2]\,.
$$
\item[(b)] For a differentiable function  $g:\mathbb{R}^+\mapsto\mathbb{R}$ such that $\mathrm{Var}[g(W)]$ exists, assume that either 
\begin{enumerate}
\item[(i)] $W$ is NBUE and $g(x)+xg^\prime(x)$ is decreasing in $x$; or
\item[(ii)] $W$ is NWUE and $g(x)+xg^\prime(x)$ is increasing in $x$.
\end{enumerate}
Then
$$
\mathrm{Var}[g(W)]\geq\frac{(\mathbb{E}[Wg^\prime(W)])^2}{\lambda^2\mathrm{Var}[W]}\,.
$$
\end{enumerate}
\end{theorem}
\begin{example}
Consider the random sum $W=\sum_{i=1}^NX_i$, where $X,X_1,X_2,\ldots$ are independent and identically distributed, continuous, real-valued random variables and $N$ is a counting random variable supported on the non-negative integers. Conditions are known under which $W$ is NWUE. For example, \cite{Brown90} shows that if $N$ is Geometric, then $W$ is NWUE, regardless of the distribution of $X$. More generally, Corollary 2.1 of \cite{Willmot05} establishes that if $N$ satisfies
\begin{equation}\label{eq:randomsum}
\sum_{k=0}^\infty\mathbb{P}(N>n+k+1)\geq\mathbb{P}(N>n)\sum_{k=0}^\infty\mathbb{P}(N>k)\,,
\end{equation}
for all $n=0,1,\ldots$, then $W$ is NWUE. This includes, for example, the case where $N$ is mixed Poisson with a mixing distribution that is itself NWUE; see Corollary 3.1 of \cite{Willmot05}. Thus, under the condition (\ref{eq:randomsum}), the bounds of the NWUE cases of Theorem \ref{theorem:equilibrium} apply, with $\lambda^{-1}=\mathbb{E}[N]\mathbb{E}[X]$ and $\mathrm{Var}[W]=(\mathbb{E}[X])^2\mathrm{Var}[N]+\mathbb{E}[N]\mathrm{Var}[X]$.
\end{example}
\subsection*{Acknowledgements}
Part of this work was completed while FD and YS were attending the
Workshop on New Directions in Stein's Method, held at the Institute
for Mathematical Sciences, National University of Singapore in May
2015.  We thank the IMS, and the organisers of that workshop, for
their support and hospitality.  FD also thanks the University of
Li\`ege for supporting a visit there. The research of FG and CL is
supported by a BOF Starting Grant of Ghent University.

\appendix

\section{Example \ref{eg:smoothing}: Proofs of claims}
\label{sec:proofs}

We begin by showing that $\tau_{\epsilon}(x)$, as defined in (\ref{eq:smooth}), is the Stein kernel of $Y+Z$. To see this, note that $\mathbb{P}(Y+Z\leq t)=\mathbb{E}[\Phi_\epsilon(t-Y)]$, so that $Y+Z$ has density $p_{\epsilon}(t)=\mathbb{E}[\varphi_\epsilon(t-Y)]$. Hence, since $Y+Z$ has expectation $\mu$, its Stein kernel is given by
$$
\frac{1}{p_{\epsilon}(x)}\int_x^\infty(y-\mu)p_{\epsilon}(y)\,dy
=\frac{1}{p_{\epsilon}(x)}\int_x^\infty\int_{-\infty}^\infty(y-\mu)\varphi_\epsilon(y-t)\,dF(t)\,dy\,,
$$
where $F$ is the distribution function of $Y$; see \cite{CPU94}. Applying Fubini's theorem, this is equal to
$$
\frac{1}{p_{\epsilon}(x)}\int_{-\infty}^\infty\int_{x-t}^\infty(s+t-\mu)\varphi_\epsilon(s)\,ds\,dF(t)
=\frac{1}{p_{\epsilon}(x)}\mathbb{E}\left[\epsilon^2\varphi_\epsilon(x-Y)+(Y-\mu)\bar{\Phi}_\epsilon(x-Y)\right]\,,
$$
since $\int_y^\infty s\varphi_\epsilon(s)\,ds=\epsilon^2\varphi_\epsilon(y)$.  This Stein kernel is easily seen to be equal to $\tau_{\epsilon}(x)$ given in (\ref{eq:smooth}).

Now, to prove claim (i), we firstly note that $Y\leq_{cx}Y+Z$ (see Theorem 3.A.34 of \cite{ss07}), so that $\mathbb{E}[\phi(Y)]\leq\mathbb{E}[\phi(Y+Z)]$ for any convex function $\phi$.  Noting that the function $f(\alpha)=\mathbb{E}[(g(Y)-\alpha)^2]$ is minimized at $\alpha=\mathbb{E}[g(Y)]$, we have   
$$
\mathrm{Var}[g(Y)]=\mathbb{E}\left[\left(g(Y)-\mathbb{E}[g(Y)]\right)^2\right]\leq
\mathbb{E}\left[\left(g(Y)-\mathbb{E}[g(Y+Z)]\right)^2\right]\leq\mathrm{Var}[g(Y+Z)]\,,
$$
where the final inequality follows from the assumption in (i) that the mapping $x\mapsto\left(g(x)-\mathbb{E}[g(Y+Z)]\right)^2$ is convex. Applying the upper bound from (\ref{eq:8}) completes the proof of (i).

We use a similar argument for (ii). We have that
$$
\mathrm{Var}[g(Y+Z)]\leq\mathbb{E}[(g(Y+Z)-\mathbb{E}[g(Y)])^2]\leq\mathbb{E}[(g(Y)-\mathbb{E}[g(Y)])^2]\,,
$$
where the final inequality uses the convex ordering between $Y$ and $Y+Z$ (from which $\mathbb{E}[\phi(Y+Z)]\leq\mathbb{E}[\phi(Y)]$ for any concave function $\phi$) and the assumption that the mapping $x\mapsto (g(x)-\mathbb{E}[g(Y)])^2$ is concave. We now apply the lower bound from (\ref{eq:8}) to complete the proof of (ii).

\section{Example \ref{eg:Pearson}: Stein kernel and further applications}
\label{sec:more-examples-1}
We start by recalling a result taken from \cite[Equation (40), p.\
65]{Stein1986}, which was used in Example \ref{eg:Pearson}.
 
\begin{proposition}[Pearson distribution]\label{prop:perason}
A random variable with mean $\mu$ and variance $\sigma^2$ is of
Pearson type if and only if there exist
$\delta_1, \delta_2, \delta_3 \in \R$, not all equal to 0, such that
\begin{equation*}
  \frac{p'(x)}{p(x) } = -\frac{(2\delta_1+1)(x-\mu) +
   \delta_2}{\delta_1(x-\mu)^2+\delta_2(x-\mu)+\delta_3}. 
\end{equation*}
In this case, its Stein kernel is
$ \tau(x) = \delta_1 (x-\mu)^2 + \delta_2 (x-\mu) + \delta_3.  $
\end{proposition}

To complement Example \ref{eg:Pearson} and illustrate the scope of its application, we use the remainder of this appendix to present further examples along similar lines. 

\begin{example}[Negative binomial data, inference on proportion, Beta  prior]
  If $X \sim {NB}(r, \theta)$ has a negative binomial distribution with $\theta\in[0,1]$ and fixed
  $r\in\N$, and $\Theta_0 \sim{Beta}(\alpha, \beta)$ with
  $\alpha,\beta>0$, then
  $\Theta_2 \sim
  {Beta}\left(\sum_{i=1}^nx_i+\alpha,nr+\beta\right)$. The Stein
  kernel for this Beta distribution is
  $\tau (\theta) = \frac{\theta (1 - \theta)}{\sum_{i=1}^nx_i+nr+
    \alpha + \beta}$. Consequently,
$$
\frac{(\sum_{i=1}^nx_i+nr+\alpha +\beta+1)}{(\sum_{i=1}^nx_i+\alpha)(nr + \beta)} \mathbb{E} [\Theta_2 (1-\Theta_2) g'(\Theta_2)]^2 \leq \mathrm{Var}[g (\Theta_2)] \leq \frac{\mathbb{E}[\Theta_2 (1 - \Theta_2)g'(\Theta_2)^2]}{\sum_{i=1}^nx_i+nr+\alpha +\beta}.
$$

\end{example}

\begin{example}[Weibull data, inference on scale, Inverse Gamma  prior]
  If $X \sim {Wei}(k, \theta)$ has a Weibull distribution with $\theta>0$ and fixed $k>0$   	(note that here we consider the Weibull density
  $x\mapsto \frac{k x^{k-1}}{\theta}\exp(-x^k/\theta), x>0$), and
  $\Theta_0 \sim{IG}(\alpha, \beta)$ with $\alpha,\beta>0$, then
  $\Theta_2 \sim
  {IG}\left(n+\alpha,\sum_{i=1}^nx_i^k+\beta\right)$. The Stein kernel
  for this Inverse Gamma distribution is
  $\tau (\theta) = \frac{\theta^2}{n+\alpha-1}$. Consequently,
$$
 \frac{n + \alpha -2}{(\sum_{i=1}^n x_i^k + \beta)^2} \mathbb{E}
 [\Theta_2^2 g'(\Theta_2)]^2 \leq \mathrm{Var}[g(\Theta_2)] \leq \frac{\mathbb{E} [\Theta_2^2 g'(\Theta_2)^2]}{n + \alpha -1}.
$$
\end{example}

\begin{example}[Gamma data, inference on scale,  Gamma  prior]
  If $X \sim {Gam}(k, \theta)$ has a Gamma distribution with $\theta, k>0$, and
  $\Theta_0 \sim{Gam}(\alpha, \beta)$ with $\alpha,\beta>0$, then
  $\Theta_2 \sim
  {Gam}\left(nk+\alpha,\sum_{i=1}^nx_i+\beta\right)$. The Stein kernel
  for this Gamma distribution is
  $\tau (\theta) = \frac{\theta}{\sum_{i=1}^n x_i +
    \beta}$. Consequently, 
$$
\frac{\mathbb{E} [\Theta_2 g'(\Theta_2)]^2}{nk + \alpha} \leq \mathrm{Var}[g(\Theta_2)] \leq \frac{1}{\sum_{i=1}^n x_i + \beta} \mathbb{E} [\Theta_2 g'(\Theta_2)^2].
$$
\end{example}

\begin{example}[Laplace data, inference on scale, inverse gamma  prior]
  If $X \sim {Lap}(\mu, \theta)$ has a Laplace distribution with $\theta>0$ and fixed $\mu\in\R$,
  and $\Theta_0 \sim{IG}(\alpha, \beta)$ with $\alpha,\beta>0$, then
  $\Theta_2 \sim
  {IG}\left(n+\alpha,\sum_{i=1}^n|x_i-\mu|+\beta\right)$.  The
Stein kernel can readily be deduced from previous
  examples, and we get 
$$
\frac{n + \alpha -2}{(\sum_{i=1}^n |x_i - \mu| + \beta)^2} \mathbb{E} \left[ \Theta_2^2 g'(\Theta_2) \right]^2 \leq  \mathrm{Var}[g(\Theta_2)]  \leq \frac{1}{n + \alpha -1} \mathbb{E} \left[\Theta_2^2 g'(\Theta_2)^2 \right].
$$
\end{example}

\begin{example}[Poisson data, inference on mean=scale,  Gamma  prior]
If  $X \sim {Poi}(\theta)$ has a Poisson distribution with $\theta>0$, and $\Theta_0 \sim{Gam}(\alpha, \beta)$ with $\alpha,\beta>0$, then
$\Theta_2 \sim {Gam}\left(\sum_{i=1}^nx_i+\alpha,n+\beta\right)$. The
Stein kernel can readily be
deduced from previous examples,  and we get 
$$
\frac{\mathbb{E} [\Theta_2 g'(\Theta_2)]^2}{\sum_{i=1}^nx_i+\alpha} \leq \mathrm{Var}[g(\Theta_2)] \leq  \frac{1}{n+\beta} \mathbb{E}\left[  \Theta_2 g'(\Theta_2)^2   \right].
$$
\end{example}

\begin{example}[Uniform data, inference on interval length,  Pareto  prior]
  If $X \sim {U}(0,\theta)$ has a Uniform distribution with $\theta>0$, and
  $\Theta_0 \sim{Par}(\alpha, \beta)$ has a Pareto distribution with $\alpha,\beta>0$ (as a
  reminder, the density of such a Pareto distribution is
  $\theta\mapsto
  \frac{\alpha\beta^\alpha}{\theta^{\alpha+1}}\mathbb{I}[\beta\leq\theta]$
  where $\mathbb{I}[A]$ is the indicator function of the event
  $A$), then
  $\Theta_2 \sim {Par}\left(n+\alpha,\max(m(x),\beta)\right)$ with
  $m(x)=\max(x_1,\ldots,x_n)$. The Stein kernel for this Pareto
  distribution is
  $\tau (\theta) = \frac{\max(m(x) , \beta) - \theta}{n+\alpha - 1}
  \theta$. Consequently, we get 
  \begin{multline*}
    \frac{(n + \alpha - 2)}{(n + \alpha)(\max(m(x) , \beta))^2}
    \mathbb{E} \left[(\max(m(x) , \beta)-\Theta_2)\Theta_2 g'(\Theta_2)
    \right]^2 \\
     \leq \mathrm{Var}[g(\Theta_2)]
     \leq  \frac{1}{n+\alpha - 1}\mathbb{E}\left[ (\max(m(x) , \beta)-\Theta_2) \Theta_2 g'(\Theta_2)^2 \right].
  \end{multline*}
\end{example}

\

\end{document}